\documentclass{ieja}
\usepackage{amsmath,amsthm,amssymb,amsfonts,amscd, array}
\usepackage[frame,ps,matrix,arrow,curve,rotate,all,2cell,tips]{xy}
\def\currentvolume{*}
\def\currentyear{2010}
\def\pagerange{**-**}

\newcommand{\fhead}{\hat{f}}
\newcommand{\Ubar}{\overline{U}}
\newcommand{\Ibar}{\overline{I}}

\newcommand{\catc}{\mathcal{C}}

\newcommand{\bk}{\mathbf{k}}
\newcommand{\bM}{\mathbf{M}(2)}
\newcommand{\bA}{\mathbf{A}_2}

\newcommand{\Mod}{\mathbf{Mod}}
\newcommand{\Hommod}{\mathbf{HomMod}}
\newcommand{\HNA}{\mathbf{HNA}}
\newcommand{\HA}{\mathbf{HA}}
\newcommand{\HL}{\mathbf{HL}}

\newcommand{\relphantom[1]}{\mathrel{\phantom{#1}}}

\newtheorem{theorem}{Theorem}[section]
 \newtheorem{corollary}[theorem]{Corollary}
 \newtheorem{lemma}[theorem]{Lemma}
 \newtheorem{proposition}[theorem]{Proposition}
\theoremstyle{definition}
\newtheorem{definition}[theorem]{Definition}
\newtheorem{example}[theorem]{Example}
\newtheorem{remark}[theorem]{Remark}

\title{Hom-bialgebras and comodule Hom-algebras}
\author{Donald Yau}
\address{
Department of Mathematics \\
The Ohio State University at Newark \\
1179 University Drive \\
Newark, OH 43055, USA \\
e-mail: dyau@math.ohio-state.edu}

\begin{document}

\begin{abstract}
We study Hom-bialgebras and objects admitting coactions by Hom-bialgebras.  In particular, we construct a Hom-bialgebra $\bM$ representing the functor of $2\times 2$-matrices on Hom-associative algebras.  Then we construct a Hom-algebra analogue of the affine plane and show that it is a comodule Hom-algebra over $\bM$ in a suitable sense.  It is also shown that the enveloping Hom-associative algebra of a Hom-Lie algebra is naturally a Hom-bialgebra.
\\[+2mm]
\subjclassname{ 16S30, 16W30, 17A30.}\\
{\bf Keywords}: Hom-bialgebras, comodule Hom-algebras, affine plane.
\end{abstract}

\maketitle

\section{Introduction}

A \emph{Hom-Lie algebra} $L$ has a skew-symmetric bracket $[-,-] \colon L \otimes L \to L$ and a linear map $\alpha \colon L \to L$ such that an $\alpha$-twisted version of the Jacobi identity holds.  It is said to be \emph{multiplicative} if, in addition, $\alpha \circ [-,-] = [-,-]\circ \alpha^{\otimes 2}$.  An ordinary Lie algebra is a multiplicative Hom-Lie algebra with $\alpha = Id$.  Hom-Lie algebras were introduced in \cite{hls} to describe the structures on certain deformations of the Witt algebra and the Virasoro algebra.  Earlier precursors of Hom-Lie algebras can be found in \cite{as,hu,liu}.  Hom-Lie algebras are also related to deformed vector fields \cite{hls,ls,ls2,ls3,rs,ss}, number theory \cite{larsson}, the various versions of the Yang-Baxter equations, braid group representations, and quantum groups \cite{yau5} - \cite{yau12}.  Other papers concerning Hom-type structures are \cite{am}- \cite{gohr}, \cite{jl}, \cite{makhlouf} - \cite{ms4}, and \cite{yau} - \cite{yau4}.

A fundamental property of Lie algebras is that they are related to associative algebras via the commutator bracket construction and the universal enveloping algebra functor.  The associative type objects corresponding to Hom-Lie algebras, called \emph{Hom-associative algebras}, were introduced in \cite{ms}.  A Hom-associative algebra $A$ has a multiplication $\mu \colon A \otimes A \to A$ and a linear self-map $\alpha$ such that the following $\alpha$-twisted version of associativity holds: $\mu(\alpha(x),\mu(y,z)) = \mu(\mu(x,y),\alpha(z))$.  It is said to be \emph{multiplicative} if, in addition, $\alpha \circ \mu = \mu \circ \alpha^{\otimes 2}$.  Ordinary associative algebras are multiplicative Hom-associative algebras with $\alpha = Id$.  It is proved in \cite{ms} that a Hom-associative algebra $A$ gives rise to a Hom-Lie algebra $HLie(A)$ via the commutator bracket construction, i.e., $[x,y] = \mu(x,y) - \mu(y,x)$.  Moreover, there is an enveloping Hom-associative algebra functor $U_{HLie}$ going the other way, from Hom-Lie algebras to Hom-associative algebras \cite{yau}.  As in the classical case, the functor $U_{HLie}$ is the left adjoint of the functor $HLie$.

The purpose of this paper is to advance the theory of bialgebras in the Hom-algebra setting.  Hom-type analogues of Hopf algebras have been studied in \cite{ms2}.  In particular, the authors of \cite{ms2} considered Hom-versions of the convolution product and primitive elements.  Moreover, in \cite{cg} Hom-Hopf algebras are studied from a categorical view-point.  Hom-type analogues of quantum groups, Lie bialgebras, and infinitesimal bialgebras are studied in \cite{yau7} - \cite{yau11}.

In the rest of this introduction, let us describe some of the results of this paper.  In the theory of bialgebras and quantum groups, a fundamental object is the bialgebra $M(2)$, which represents the functor $M_2(-)$ sending an associative algebra $A$ to the algebra $M_2(A)$ of $2 \times 2$-matrices with elements in $A$.  Moreover, the matrix multiplication $M_2(A) \times M_2(A) \to M_2(A)$ is universally represented by the comultiplication on $M(2)$.

We will define a suitable notion of \emph{Hom-bialgebra}, in which the comultiplication $\Delta$ satisfies an $\alpha$-twisted version of coassociativity and is a morphism of Hom-associative algebras.  In the first main result of this paper, we construct a Hom-bialgebra $\bM$, which is the analogue of $M(2)$ in the Hom-algebra setting.  We show that $\bM$ is the representing object of the functor that sends a multiplicative Hom-associative algebra $A$ to the multiplicative Hom-associative algebra $M_2(A)$ of $2 \times 2$-matrices with elements in $A$.  Moreover, the multiplication on $M_2(A)$ is universally represented by the comultiplication on $\bM$.

A more precise version of the following result is proved in Sections ~\ref{sec:homas} and ~\ref{sec:hombi}.

\begin{theorem}
\label{thm:M2}
There exists a Hom-bialgebra $\bM$ such that:
\begin{enumerate}
\item
For any multiplicative Hom-associative algebra $A$, there is a natural bijection
\[
\HA(\bM,A) \cong M_2(A),
\]
where $\HA$ is the category of multiplicative Hom-associative algebras and $\HA(-,-)$ denotes the morphism sets in $\HA$.
\item
The matrix multiplication on $M_2(A)$ is represented by the comultiplication on $\bM$.
\end{enumerate}
\end{theorem}

The significance of the bialgebra $M(2)$ goes beyond the functor $M_2$ represented by it.  In fact, the affine plane $A_2 = \bk[x,y]$ is a non-trivial $M(2)$-comodule algebra.  In other words, there is a non-trivial $M(2)$-comodule structure
\[
\rho \colon \bk[x,y] \to M(2) \otimes \bk[x,y]
\]
on the affine plane such that $\rho$ is a morphism of algebras. We will construct a multiplicative Hom-associative algebra $\bA$, the \emph{Hom-affine plane}, which is the Hom-algebra analogue of the affine plane.  We define what it means for a multiplicative Hom-associative algebra to be a \emph{comodule Hom-algebra} over a Hom-bialgebra.

A more precise version of the following result is proved in Section ~\ref{sec:comod}.

\begin{theorem}
\label{thm:comod}
The Hom-affine plane $\bA$ admits the structure of a non-trivial $\bM$-comodule Hom-algebra.
\end{theorem}

The rest of this paper is organized as follows.  In Section ~\ref{sec:homas} we first recall some basic definitions about Hom-algebras.  The free functor from modules to multiplicative Hom-associative algebras is constructed (Corollary ~\ref{cor:F}).  The multiplicative Hom-associative algebra $\bM$ is constructed using this free functor.  It is then shown that $\bM$ represents the functor $M_2$ of $2 \times 2$-matrices on multiplicative Hom-associative algebras (Corollary ~\ref{cor1:M2}).

In Section ~\ref{sec:hombi}, we equip $\bM$ with a comultiplication $\Delta$ and observe that $\Delta$ represents the matrix multiplication on $M_2(A)$ for a multiplicative Hom-associative algebra $A$.  Then we observe that, with its comultiplication, $\bM$ forms a Hom-bialgebra (Corollary ~\ref{M2bialg}).  It is also observed that a unital version of the enveloping Hom-associative algebra $U_{HLie}(L)$ is a Hom-bialgebra (Theorem ~\ref{U}).

In Section ~\ref{sec:comod}, we define the Hom-affine plane $\bA$, using once again the free functor in Corollary ~\ref{cor:F}.  Then we extend the usual notion of comodule algebra to the Hom-algebra setting.  Finally, we show that the Hom-affine plane $\bA$ is a non-trivial $\bM$-comodule Hom-algebra (Theorem ~\ref{thm:comodule}).

\section{The Hom-associative algebra $\bM$}
\label{sec:homas}

Before we construct $\bM$, we first recall some relevant definitions about Hom-modules and Hom-nonassociative algebras that were first introduced in \cite{yau}.  The multiplicative Hom-associative algebra $\bM$ will be constructed using the free multiplicative Hom-associative algebra functor $F$ in Corollary ~\ref{cor:F}.  The functor $F$ is the composition of several free functors, which we discuss first.

\subsection{Conventions}
Throughout the rest of this paper, let $\bk$ denote a field of characteristic $0$.  Unless otherwise specified, linearity, modules, and $\otimes$ are all meant over $\bk$.  The category of $\bk$-modules is denoted by $\Mod$.  The free module on a set $S$ is denoted by $\bk\langle S\rangle$.  If $S = \{a_1, \ldots , a_n\}$, then we write $\bk\langle S\rangle$ as $\bk\langle a_1, \ldots , a_n\rangle$.  Given two $\bk$-modules $M$ and $N$, denote by $\tau \colon M \otimes N \cong N \otimes M$ the twist isomorphism, i.e., $\tau(m \otimes n) = n \otimes m$ for $m \in M$ and $n \in $N.

For a category $\catc$ and two objects $x$ and $y$ in $\catc$, we denote by $\catc(x,y)$ the set of morphisms from $x$ to $y$ in $\catc$.

\subsection{Hom-modules}
\label{subsec:Hom-mod}

By a \textbf{Hom-module}, we mean a pair $(V, \alpha)$ in which $V$ is a module and $\alpha \colon V \to V$ is a linear map.  A \textbf{morphism} $(V, \alpha) \to (V^\prime, \alpha^\prime)$ of Hom-modules is a linear map $f \colon V \to V^\prime$ such that $\alpha^\prime \circ f = f \circ \alpha$.  The category of Hom-modules is denoted by $\Hommod$.  When there is no danger of confusion, we will denote a Hom-module $(V,\alpha)$ simply by $V$.

There is a forgetful functor $U \colon \Hommod \to \Mod$ that sends a Hom-module $(V,\alpha)$ to the module $V$, forgetting about the map $\alpha$.  Conversely, there is a free Hom-module associated to a module.

\begin{theorem}
\label{thm:F0}
There is an adjoint pair of functors
\[
F_0 \colon \Mod \rightleftarrows \colon \Hommod \colon U
\]
in which $F_0$ is the left adjoint.
\end{theorem}

\begin{proof}
Let $V = \oplus \bk\langle S\rangle$ be a module.  For each $x \in S$, let $x_i$ $(i \geq 1)$ be a sequence of independent variables.  We also set $x = x_0$.  Then we define the module
\[
F_0(V) = \bigoplus_{x \in S;\, i \geq 0} \bk\langle x_i\rangle.
\]
Define a linear map $\alpha \colon F_0(V) \to F_0(V)$ by setting $\alpha(x_i) = x_{i+1}$ for each $x \in S$ and $i \geq 0$.  Then $(F_0(V),\alpha)$ is a Hom-module.  Let $\iota \colon V \to F_0(V)$ be the inclusion map defined as $\iota(x) = x_0$ for $x \in S$.

Let $(N,\alpha_N)$ be a Hom-module, and let $f \colon V \to N$ be a linear map.  We must show that there exists a unique morphism $\fhead \colon F_0(V) \to N$ of Hom-modules such that $f = \fhead \circ \iota$.  The desired map $\fhead$ is defined on the generators as
\[
\fhead(x_i) = \alpha_N^i(f(x))
\]
for $x \in S$ and $i \geq 0$, where $\alpha_N^i = \alpha_N \circ \cdots \circ \alpha_N$ with $i$ copies of $\alpha_N$.  It is clear from this definition that $f = \fhead \circ \iota$.  It is compatible with $\alpha$ because
\[
\alpha_N(\fhead(x_i)) = \alpha_N^{i+1}(f(x)) = \fhead(x_{i+1}) = \fhead(\alpha(x_i)).
\]
This map $\fhead$ is unique.  Indeed, let $g \colon F_0(V) \to N$ be another morphism of Hom-modules that satisfies $f = g \circ \iota$.  Then we have
\[
g(x_i) = g(\alpha^i(x_0)) = \alpha_N^i(g(x_0)) = \alpha_N^i(f(x)) = \fhead(x_i).
\]
This finishes the proof.
\end{proof}


\subsection{Hom-nonassociative algebras}
\label{subsec:HomNA}

By a \textbf{Hom-nonassociative algebra}, we mean a quadruple $(A, \mu, \alpha,\eta)$ in which:
\begin{enumerate}
\item $(A,\alpha)$ is a Hom-module;
\item $\mu \colon A \otimes A \to A$, the multiplication, is a bilinear map;
\item
$\eta \colon \bk \to A$, the unit, is a linear map such that the following diagram commutes:
\[
\SelectTips{cm}{10}
\xymatrix{
\bk \otimes A \ar[r]^-{\eta \otimes Id} \ar[dr]_-{\cong} & A \otimes A \ar[d]^-{\mu} & A \otimes \bk \ar[l]_-{Id \otimes \eta} \ar[dl]^-{\cong}\\
 & A &
}
\]
\end{enumerate}
We will write $\eta(1) \in A$ as $1_A$ or simply $1$.  We also abbreviate $\mu(x,y)$ to $xy$ if there is no danger of confusion.

A morphism $f \colon (A, \mu, \alpha,\eta) \to (A^\prime, \mu^\prime, \alpha^\prime,\eta')$ of Hom-nonassociative algebras is a morphism $f \colon (A,\alpha) \to (A^\prime,\alpha')$ of Hom-modules such that $f \circ \mu = \mu^\prime \circ f^{\otimes 2}$ and $f\circ\eta = \eta'$.  The category of Hom-nonassociative algebras is denoted by $\HNA$.

\begin{remark}
\label{remark:nonunital}
Note that in \cite{yau}, a Hom-nonassociative algebra is not required to have a unit $\eta$.  Here we call them \textbf{non-unital Hom-nonassociative algebras}.  A morphism $f \colon (A,\mu,\alpha) \to (A',\mu',\alpha')$ of non-unital Hom-nonassociative algebras is defined as in the previous paragraph, omitting the condition $f\circ\eta = \eta'$.  The category $\HNA$ is a faithful subcategory of the category of non-unital Hom-nonassociative algebras.  In what follows, it should be clear from the context whether we consider units or not.
\end{remark}

Suppose that $(A, \mu, \alpha,\eta)$ and $(A', \mu', \alpha',\eta')$ are Hom-nonassociative algebras.  Then their tensor product is defined in the usual way, $(A \otimes A', \mu'', \alpha \otimes \alpha', \eta \otimes \eta')$.  Here $\mu''$ is the composition
\[
(A \otimes A') \otimes (A \otimes A') \xrightarrow{Id_A \otimes \tau \otimes Id_{A'}} A^{\otimes 2} \otimes A'^{\otimes 2} \xrightarrow{\mu \otimes \mu'} A \otimes A',
\]
where $\tau$ is the twist isomorphism $A' \otimes A \cong A \otimes A'$.

There is a forgetful functor $U \colon \HNA \to \Hommod$ that forgets about the multiplication $\mu$ and the unit $\eta$.  Conversely, there is a free Hom-nonassociative algebra associated to each Hom-module.

\begin{theorem}
\label{thm:F1}
There is an adjoint pair of functor
\[
F_1 \colon \Hommod \rightleftarrows \HNA \colon U
\]
in which $F_1$ is the left adjoint.
\end{theorem}

\begin{proof}
This is, in fact, just \cite[Theorem 1, p.100]{yau} with very minor modifications to account for the unit $\eta$.  Using the notations in \cite{yau}, if $(V,\alpha)$ is a Hom-module, then we define the module
\[
F_1(V) = \bk \oplus F_{HNAs}(V),
\]
where $F_{HNAs}$ is the free non-unital Hom-nonassociative algebra functor in \cite[Theorem 1]{yau}.  In more details, it is defined as
\begin{equation}
\label{FHNAs}
F_{HNAs}(V) = \bigoplus_{n \geq 1} \bigoplus_{\tau \in T^{wt}_n} V_{\tau}^{\otimes n},
\end{equation}
where $T^{wt}_n$ is the set of \emph{weighted $n$-trees} \cite[section 2.3]{yau} and $V_{\tau}^{\otimes n}$ is a copy of $V^{\otimes n}$.  Since there is only one weighted $1$-tree, there is a natural inclusion $V \hookrightarrow F_{HNAs}(V)$.  This gives a natural inclusion $V \hookrightarrow F_1(V)$.

We extend the map $\alpha$ on $F_{HNAs}(V)$ to $F_1(V)$ by setting $\alpha_F|\bk = Id_{\bk}$.  The unit $\eta_F \colon \bk \to F_1(V)$ is the obvious inclusion.  The multiplication $\mu$ on $F_{HNAs}(V)$ is extended to $F_1(V)$ by setting
\[
\mu_F((a,x);(b,y)) = (ab, ay + bx + \mu(x,y))
\]
for $a, b \in \bk$ and $x, y \in F_{HNAs}(V)$.  It is straightforward to check that
\[
(F_1(V),\mu_F,\alpha_F,\eta_F)
\]
is a Hom-nonassociative algebra.  The proof that $F_1$ is the left adjoint of the forgetful functor is essentially identical to the proof of \cite[Theorem 1]{yau}.
\end{proof}


\subsection{Hom-associative algebra}
\label{subsec:HA}


By a \textbf{Hom-associative algebra}, we mean a Hom-nonassociative algebra $(A,\mu,\alpha,\eta)$ such that
\[
\mu \circ (\mu \otimes \alpha) = \mu \circ (\alpha \otimes \mu) \quad \text{(Hom-associativity)}
\]
holds.  A Hom-associative algebra is said to be \textbf{multiplicative} if, in addition, $\mu \circ \alpha^{\otimes 2} = \alpha \circ \mu$ (multiplicativity).  A morphism of (multiplicative) Hom-associative algebras is a morphism of the underlying Hom-nonassociative algebras.  The category of multiplicative Hom-associative algebras is denoted by $\HA$, which is a full and faithful subcategory of the category $\HNA$ of Hom-nonassociative algebras.  As before, \textbf{non-unital (multiplicative) Hom-associative algebras} are non-unital Hom-nonassociative algebras satisfying Hom-associativity (and multiplicativity).

There is a subcategory inclusion $U \colon \HA \to \HNA$.  Conversely, every Hom-nonassociative algebra has a unique largest multiplicative Hom-associative algebra quotient.

\begin{theorem}
\label{thm:F2}
There is an adjoint pair of functors
\[
F_2 \colon \HNA \rightleftarrows \HA \colon U
\]
in which $F_2$ is the left adjoint.
\end{theorem}

\begin{proof}
This is, in fact, \cite[section 4.3, p.103]{yau} with very minor modifications to account for the multiplicativity of $\alpha$.  Let $(A,\mu,\alpha,\eta)$ be a Hom-nonassociative algebra.  Using the notations in \cite[section 4.3]{yau}, we define inductively:
\[
\begin{split}
J^1 &= \langle im(\mu\circ(\mu\otimes\alpha - \alpha\otimes\mu)); im(\alpha \circ \mu - \mu\circ\alpha^{\otimes 2}) \rangle,\\
J^{n+1} &= \langle J^n \cup \alpha(J^n)\rangle, \quad \text{and} \quad J^{\infty} = \bigcup_{n \geq 1}J^n.
\end{split}
\]
For a subset $S \subseteq A$, the notation $\langle S \rangle$ denotes the smallest submodule of $A$ containing $S$ such that $\mu(\langle S\rangle,A) \subseteq \langle S\rangle$ and $\mu(A,\langle S\rangle) \subseteq \langle S\rangle$ \cite[section 3.4]{yau}.  Now we define the quotient module
\[
F_2(A) = A/J^\infty,
\]
equipped with the induced maps of $\mu$ and $\alpha$.  The unit $\eta_F \colon \bk \to F_2(A)$ of $F_2(A)$ is the composition $\bk \to A \to F_2(A)$ of the unit $\eta$ of $A$ and the quotient map.  It is easy to check that $F_2(A)$ is a Hom-associative algebra and that $F_2$ is the left adjoint of the subcategory inclusion $U$.
\end{proof}

Note that $F_2 \circ U = Id_{\HA}$, since for a multiplicative Hom-associative algebra $A$, all the $J^n = 0$.  Thus, Theorem ~\ref{thm:F2} says that $\HA$ is a reflective subcategory of $\HNA$.

From Theorems ~\ref{thm:F0}, \ref{thm:F1}, and \ref{thm:F2}, we have the following three adjoint pairs of functors:
\begin{equation}
\label{adjoints}
\SelectTips{cm}{10}
\xymatrix{
\Mod\ar@<.5ex>[r]^-{F_0} & \Hommod \ar@<.5ex>[r]^-{F_1} \ar@<.5ex>[l] & \HNA \ar@<.5ex>[r]^-{F_2} \ar@<.5ex>[l] & \HA \ar@<.5ex>[l].
}
\end{equation}
The three $F_i$ are the left adjoints, and the right adjoints are all forgetful functors.  Since adjoint pairs can be composed to yield another adjoint pair, we obtain the following immediate consequence.

\begin{corollary}
\label{cor:F}
There is an adjoint pair of functors
\[
F = F_2 \circ F_1 \circ F_0 \colon \Mod \rightleftarrows \HA \colon U,
\]
in which $F$ is the left adjoint and $U(A,\mu,\alpha,\eta) = A$.
\end{corollary}

We call $F \colon \Mod \to \HA$ the \textbf{free multiplicative Hom-associative algebra functor}.

For a module $V$, using the adjoint pairs in \eqref{adjoints}, there are natural maps
\begin{equation}
\label{j}
V \hookrightarrow F_0(V) \hookrightarrow F_1(F_0(V)) \twoheadrightarrow F_2(F_1(F_0(V))) = F(V).
\end{equation}
Denote by $j_V \colon V \to F(V)$ the composition of these maps.  For an element $a \in V$, we will abbreviate $j_V(a) \in F(V)$ to $a$.

\subsection{Matrix Hom-associative algebras}

For a module $V$, let $M_2(V)$ denote the module of $2 \times 2$-matrices with elements in $V$.

\begin{lemma}
\label{lem:M2A}
Let $(A,\mu,\alpha,\eta)$ be a multiplicative Hom-associative algebra.  Then $(M_2(A),\mu',\alpha',\eta')$ is also a multiplicative Hom-associative algebra in which: \begin{enumerate}
\item
the multiplication $\mu'$ is given by matrix multiplication and $\mu$;
\item
$\alpha'$ is given by $\alpha$ in each entry;
\item
$\eta'(1) = \begin{pmatrix} 1_A & 0 \\ 0 & 1_A\end{pmatrix}$.
\end{enumerate}
\end{lemma}

\begin{proof}
To see that $\mu'$ satisfies Hom-associativity, let $X = (x_{ij})$, $Y = (y_{ij})$, and $Z = (z_{ij})$ be elements in $M_2(A)$.  Then the $(i,j)$-entries in $(XY)\alpha'(Z)$ and $\alpha'(X)(YZ)$ are
\[
\sum_{1 \leq l,k \leq 2} (x_{ik}y_{kl})\alpha(z_{lj}) \quad\text{and}\quad
\sum_{1 \leq l,k \leq 2} \alpha(x_{ik})(y_{kl}z_{lj}),
\]
respectively.  Since $A$ is a Hom-associative algebra, we have $(xy)\alpha(z) = \alpha(x)(yz)$, which shows that $\mu'$ satisfies Hom-associativity.  The other axioms for $M_2(A)$ to be a multiplicative Hom-associative algebra are equally easy to check.
\end{proof}

\begin{definition}
\label{def:M2}
Let $a$, $b$, $c$, and $d$ be independent variables.  Define
\[
\bM = F(\bk\langle a,b,c,d\rangle),
\]
where $F$ is the free multiplicative Hom-associative algebra functor in Corollary ~\ref{cor:F}.
\end{definition}
In other words, $\bM$ is the free multiplicative Hom-associative algebra on four variables.  As in \eqref{j}, we have a natural map $j \colon \bk\langle a,b,c,d\rangle \to \bM$.  With a slight abuse of notation, we will use $a$, $b$, $c$, and $d$ to denote their images under $j$ as well.

We will use the following result to deduce that $\bM$ is the representing object of the matrix functor $M_2$ on the category $\HA$ of multiplicative Hom-associative algebras.

\begin{theorem}
\label{main1}
Let $S$ be a set, and let $A$ be a multiplicative Hom-associative algebra.  Then there is a natural bijection
\[
\HA(F(\bk\langle S \rangle),A) \cong \prod_{x\in S} A
\]
given by $f \mapsto (f(x))_{x \in S}$ for $f \in \HA(F(\bk\langle S \rangle),A)$.
\end{theorem}

\begin{proof}
By the adjunction in Corollary ~\ref{cor:F}, we have natural bijections
\[
\HA(F(\bk\langle S \rangle),A)  \cong \Mod(\bk\langle S\rangle,A) \cong \prod_{x \in S}  \Mod(\bk\langle x\rangle,A) \cong  \prod_{x \in S} A,
\]
as desired.
\end{proof}

\begin{corollary}
\label{cor1:M2}
Let $A$ be a multiplicative Hom-associative algebra.  Then there is a natural bijection
\[
\HA(\bM,A) \cong M_2(A)
\]
given by
\[
f \longmapsto \begin{pmatrix} f(a) & f(b) \\ f(c) & f(d)\end{pmatrix}
\]
for $f \in \HA(\bM,A)$.
\end{corollary}

\begin{proof}
As a set, we have $M_2(A) = A^4$.  The assertion now follows from Theorem ~\ref{main1}, since $\bM = F(\bk\langle a,b,c,d\rangle)$.
\end{proof}

\section{The Hom-bialgebra $\bM$}
\label{sec:hombi}

The main purpose of this section is to observe that $\bM$ is not just a multiplicative Hom-associative algebra, but a \emph{Hom-bialgebra}.  In view of Corollary ~\ref{cor1:M2}, it is natural to ask if the matrix multiplication $\mu' \colon M_2(A) \times M_2(A) \to M_2(A)$ (Lemma ~\ref{lem:M2A}) can be universally represented by a comultiplication $\Delta$ on the representing object $\bM$.  This is true, as we will observe below (Theorem ~\ref{thm:comult}).  The comultiplication $\Delta$ makes $\bM$ into a Hom-bialgebra (Corollary ~\ref{M2bialg}).  We also observe that a unital version of the enveloping Hom-associative algebra $U_{HLie}(L)$ is a Hom-bialgebra (Theorem ~\ref{U}).

First we need some preliminary observations.

\begin{lemma}
\label{lem1:F}
Let $M$ and $N$ be two modules.  Then there is a natural isomorphism of multiplicative Hom-associative algebras,
\[
F(M) \otimes F(N) \cong F(M \oplus N),
\]
where $F$ is the free multiplicative Hom-associative algebra functor (Corollary ~\ref{cor:F}).
\end{lemma}

\begin{proof}
This is essentially identical to the proof of the same property of the universal enveloping algebra functor.  In one direction, using Corollary ~\ref{cor:F}, the morphism
\[
F(M \oplus N) \to F(M) \otimes F(N)
\]
is determined by the linear map
\[
(x,y) \mapsto x \otimes 1 + 1 \otimes y
\]
for $x \in M$ and $y \in N$.  In the other direction, the inclusions $M \hookrightarrow M \oplus N$ and $N \hookrightarrow M \oplus N$ induce morphisms $F(M) \to F(M \oplus N)$ and $F(N) \to F(M \oplus N)$, respectively, of multiplicative Hom-associative algebras.  These last two morphisms together induce the morphism $F(M) \otimes F(N) \to F(M \oplus N)$.  It is easy to check that the two morphisms that we have defined are inverse to each other.
\end{proof}

\begin{lemma}
\label{lem2:F}
Let $A$ be a multiplicative Hom-associative algebra.  Then there is a natural bijection
\[
\HA(\bM \otimes \bM,A) \cong M_2(A) \times M_2(A).
\]
\end{lemma}

\begin{proof}
Using Corollary ~\ref{cor:F}, Corollary ~\ref{cor1:M2}, and Lemma ~\ref{lem1:F}, we have the following natural bijections:
\[
\begin{split}
M_2(A) \times M_2(A)
& \cong \HA(\bM,A) \times \HA(\bM,A)\\
& \cong \Mod(\bk\langle a,b,c,d\rangle,A) \times \Mod(\bk\langle a,b,c,d\rangle,A)\\
& \cong \Mod(\bk\langle a',a'',b',b'',c',c'',d',d''\rangle,A)\\
& \cong \HA(F(\bk\langle a',a'',b',b'',c',c'',d',d''\rangle),A)\\
& \cong \HA(\bM \otimes \bM,A).
\end{split}
\]
\end{proof}

In view of the bijections in Corollary ~\ref{cor1:M2} and Lemma ~\ref{lem2:F}, it makes sense to try to represent the multiplication on $M_2(A)$ by a morphism $\Delta \colon \bM \to \bM \otimes \bM$ of multiplicative Hom-associative algebras.  Using Corollary ~\ref{cor:F} and Lemma ~\ref{lem1:F}, we have
\begin{equation}
\label{M2M2}
\begin{split}
\HA(\bM,\bM \otimes \bM)
& = \HA(F(\bk\langle a,b,c,d\rangle),\bM \otimes \bM)\\
& \cong \Mod(\bk\langle a,b,c,d\rangle,F(\bk\langle a',a'',b',b'',c',c'',d',d''\rangle)).
\end{split}
\end{equation}
So it suffices to construct a suitable linear map
\[
\Delta \colon \bk\langle a,b,c,d\rangle \to F(\bk\langle a',a'',b',b'',c',c'',d',d''\rangle).
\]
To this end, we define such a linear map $\Delta$ as
\begin{equation}
\label{Delta}
\Delta\begin{pmatrix} a & b \\ c & d\end{pmatrix} =
\begin{pmatrix}\Delta(a) & \Delta(b) \\ \Delta(c) & \Delta(d)\end{pmatrix} =
\begin{pmatrix}a' & b'\\ c' & d'\end{pmatrix}
\begin{pmatrix}a'' & b''\\ c'' & d''\end{pmatrix}.
\end{equation}
The matrix multiplication above is performed using the multiplication in the multiplicative Hom-associative algebra $F(\bk\langle a',a'',b',b'',c',c'',d',d''\rangle)$.  For example, we have $\Delta(a) = a'a'' + b'c''$ and $\Delta(c) = c'a'' + d'c''$.

\begin{theorem}
\label{thm:comult}
Let $A$ be a multiplicative Hom-associative algebra.  Under the identifications of Corollary ~\ref{cor1:M2} and Lemma ~\ref{lem2:F}, the multiplication on the multiplicative Hom-associative algebra $M_2(A)$ is the induced map
\[
\Delta^* \colon \HA(\bM \otimes \bM,A) \to \HA(\bM,A),
\]
where $\Delta \colon \bM \to \bM \otimes \bM$ is the map in \eqref{Delta}.
\end{theorem}

\begin{proof}
This is just a matter of tracing through the various bijections.
\end{proof}

The following result shows that the map $\Delta$ behaves like the dual of the multiplication in a Hom-associative algebra.

\begin{proposition}
\label{prop:Delta}
The map $\Delta \colon \bM \to \bM \otimes \bM$ in \eqref{Delta} satisfies
\[
(\Delta \otimes \alpha) \circ \Delta = (\alpha \otimes \Delta) \circ \Delta,
\]
where $\alpha \colon \bM \to \bM$ is part of the structure of the multiplicative Hom-associative algebra $\bM$.
\end{proposition}

\begin{proof}
By Corollary ~\ref{cor:F}, it suffices to check this equality on the generators $a$, $b$, $c$, and $d$.  Using the notations in \eqref{Delta}, we have
\[
\begin{split}
((\Delta \otimes \alpha) \circ \Delta)\begin{pmatrix} a & b \\ c & d\end{pmatrix}
& = \left[\begin{pmatrix} a' & b' \\ c' & d'\end{pmatrix}\begin{pmatrix} a'' & b'' \\ c'' & d''\end{pmatrix}\right]\begin{pmatrix} \alpha(a''') & \alpha(b''') \\ \alpha(c''') & \alpha(d''')\end{pmatrix}\\
& = \begin{pmatrix} \alpha(a') & \alpha(b') \\ \alpha(c') & \alpha(d')\end{pmatrix}\left[\begin{pmatrix} a'' & b'' \\ c'' & d''\end{pmatrix}\begin{pmatrix} a''' & b''' \\ c''' & d'''\end{pmatrix}\right]\\
& = ((\alpha \otimes \Delta) \circ \Delta)\begin{pmatrix} a & b \\ c & d\end{pmatrix}.
\end{split}
\]
In the second equality above, we used the Hom-associativity of the multiplication in $\bM^{\otimes 3} \cong F(\bk\langle a',a'',a''',\ldots\rangle)$ (Lemma ~\ref{lem1:F}).
\end{proof}

The previous result motivates the following definition, which is a slight modification of \cite{ms2} (Remark 3.3).

\begin{definition}
A \textbf{Hom-bialgebra} is a quadruple $(B,\mu,\Delta,\alpha)$ in which:
\begin{enumerate}
\item
$(B,\mu,\alpha)$ is a (non-unital) multiplicative Hom-associative algebra;
\item
The linear map $\Delta \colon B \to B^{\otimes 2}$, the comultiplication, is Hom-coassociative, in the sense that
\[
(\Delta \otimes \alpha) \circ \Delta = (\alpha \otimes \Delta) \circ \Delta.
\]
\item
$\Delta$ is a morphism of Hom-associative algebras.
\end{enumerate}
\end{definition}

\begin{example}
A bialgebra $(B,\mu,\Delta)$ in the usual sense is also a Hom-bialgebra with $\alpha = Id_B$.\qed
\end{example}

\begin{example}
\label{ex:bialgdeform}
More generally, let $(B,\mu,\Delta)$ be a bialgebra, and let $\alpha \colon B \to B$ be a morphism of bialgebras.  Then there is a Hom-bialgebra $B_\alpha = (B,\mu_\alpha,\Delta_\alpha,\alpha)$ in which $\mu_\alpha = \alpha\circ\mu$ and $\Delta_\alpha = \Delta \circ \alpha$.  In fact, since $B$ is a bialgebra and $\alpha$ is a bialgebra morphism, one only needs to observe that $\mu_\alpha$ and $\Delta_\alpha$ are Hom-associative and Hom-coassociative, respectively.  The Hom-associativity of $\mu_\alpha$ is shown in \cite[Corollary 2.5]{yau2}.  The Hom-coassociativity of $\Delta_\alpha$ is shown by the dual argument; see \cite[Theorem 3.16]{ms4}.  In summary, a bialgebra $B$ deforms into a Hom-bialgebra $B_\alpha$ via any bialgebra endomorphism $\alpha$.  Moreover, if $\alpha = Id$, then we have $B_{Id} = B$.\qed
\end{example}

\begin{corollary}
\label{M2bialg}
Equipped with the comultiplication $\Delta \colon \bM \to \bM \otimes \bM$ in \eqref{Delta}, the multiplicative Hom-associative algebra $\bM$ becomes a Hom-bialgebra.
\end{corollary}

\begin{proof}
This is an immediate consequence of Definition ~\ref{def:M2}, the bijections \eqref{M2M2}, and Proposition ~\ref{prop:Delta}.
\end{proof}

\subsection{Enveloping Hom-bialgebra}

The enveloping algebra $U(L)$ of a Lie algebra $L$ is a bialgebra.  In the rest of this section, we observe that a unital version of the enveloping Hom-associative algebra $U_{HLie}(L)$ of a multiplicative Hom-Lie algebra $L$ is a Hom-bialgebra.  Let us first recall some relevant definitions.

A \textbf{multiplicative Hom-Lie algebra} is a (non-unital) Hom-nonassociative algebra $(L, \lbrack -,-\rbrack, \alpha)$, satisfying the following three conditions for $x, y, z \in L$:
\begin{enumerate}
\item
$\lbrack x, y \rbrack = - \lbrack y, x \rbrack$ (skew-symmetry);
\item
$0 = \lbrack \alpha(x), \lbrack y, z \rbrack\rbrack + \lbrack \alpha(z), \lbrack x, y \rbrack \rbrack + \lbrack \alpha(y), \lbrack z, x \rbrack\rbrack$ (Hom-Jacobi identity);
\item
$\alpha([x,y]) = [\alpha(x),\alpha(y)]$ (multiplicativity of $\alpha$).
\end{enumerate}
Note that a Lie algebra is a multiplicative Hom-Lie algebra with $\alpha = Id$, since in this case the Hom-Jacobi identity reduces to the usual Jacobi identity.

A morphism of multiplicative Hom-Lie algebras is a morphism of the underlying Hom-nonassociative algebras.  The category of multiplicative Hom-Lie algebras is denoted by $\HL$.  Note that the multiplicativity of $\alpha$ is not part of the original definition of a Hom-Lie algebra in \cite{hls,ms}.

Given a (non-unital) multiplicative Hom-associative algebra $(A, \mu, \alpha)$, one can associate to it a multiplicative Hom-Lie algebra $(HLie(A), \lbrack -,-\rbrack, \alpha)$ in which $HLie(A)$ is equal to $A$ as a module and $\lbrack x, y \rbrack = xy - yx$ for $x, y \in A$ \cite[Proposition 1.7]{ms}.  This construction gives a functor $HLie$ from (non-unital) multiplicative Hom-associative algebras to multiplicative Hom-Lie algebras.  In \cite[section 4.2]{yau} the author constructed the left adjoint $U_{HLie}$ of $HLie$.  It is defined as
\[
U_{HLie}(L) = F_{HNAs}(L)/I^\infty
\]
for a multiplicative Hom-Lie algebra $L$.  Here $F_{HNAs}$ is the free non-unital Hom-nonassociative algebra functor in \eqref{FHNAs} and $I^\infty$ is a certain submodule.

We now modify $U_{HLie}$ slightly to account for the unit.  Fix a multiplicative Hom-Lie algebra $(L,[-,-],\alpha)$.  Consider the Hom-nonassociative algebra
\[
(F_1(L),\mu_F,\alpha_F,\eta_F), 
\]
obtained by applying the functor $F_1$ (Theorem ~\ref{thm:F1}) to the Hom-module $(L,\alpha)$.  As in the proof of Theorem ~\ref{thm:F2}, we consider the following sequence of submodules in $F_1(L)$:
\[
\begin{split}
\Ibar^1 & = \langle im(\mu_F \circ (\mu_F \otimes \alpha_F - \alpha_F \otimes \mu_F); im(\alpha_F \circ \mu_F - \mu_F \circ \alpha_F^{\otimes 2});\\
& \relphantom{} \relphantom{} [x,y] - (xy - yx) ~\text{for}~ x,y \in L\rangle\\
\Ibar^{n+1} & = \langle \Ibar^n \cup \alpha_F(\Ibar^n)\rangle, \quad\text{and}\quad \Ibar^\infty = \bigcup_{n \geq 1} \Ibar^n.
\end{split}
\]
It is easy to check that the module
\[
\Ubar(L) = F_1(L)/\Ibar^\infty,
\]
when equipped with the induced maps of $\mu_F$, $\alpha_F$, and $\eta_F$, is a multiplicative Hom-associative algebras.  Essentially the same proof as in \cite[section 4.2]{yau} gives the following result.

\begin{proposition}
\label{prop:Ubar}
There is an adjoint pair of functors
\[
\Ubar \colon \HL \rightleftarrows \HA \colon HLie,
\]
in which $\Ubar$ is the left adjoint.
\end{proposition}

Using the proof of Lemma ~\ref{lem1:F}, we obtain the following result.

\begin{lemma}
Let $L$ and $L'$ be multiplicative Hom-Lie algebras.  Then there is a natural isomorphism
\[
\Ubar(L) \otimes \Ubar(L') \cong \Ubar(L \oplus L')
\]
of multiplicative Hom-associative algebras.
\end{lemma}

Fix a multiplicative Hom-Lie algebra $L$.  Now we equip the multiplicative Hom-associative algebra $\Ubar(L)$ with a comultiplication.  Consider the linear map
\[
\Delta_L \colon L \to L \oplus L
\]
given by $x \mapsto (\alpha(x),\alpha(x))$.  This is a morphism of multiplicative Hom-Lie algebras.  Applying the functor $\Ubar$ above, we obtain a morphism
\begin{equation}
\label{DeltaU}
\begin{split}
\Delta_L \colon \Ubar(L)& \to \Ubar(L \oplus L) \cong \Ubar(L) \otimes \Ubar(L),\\
x& \longmapsto \alpha(x) \otimes 1 + 1 \otimes \alpha(x)
\end{split}
\end{equation}
of multiplicative Hom-associative algebras.

\begin{theorem}
\label{U}
Let $L$ be a multiplicative Hom-Lie algebra.  Equipped with the comultiplication $\Delta_L$ in \eqref{DeltaU}, the multiplicative Hom-associative algebra $\Ubar(L)$ becomes a Hom-bialgebra.
\end{theorem}

\begin{proof}
It remains to show that $\Delta_L$ is Hom-coassociative.  By Proposition ~\ref{prop:Ubar}, it suffices to check this for elements $x \in L$.  We compute as follows:
\[
\begin{split}
((\alpha \otimes \Delta_L)\circ\Delta_L)(x)
& = (\alpha \otimes \Delta_L)(\alpha(x) \otimes 1 + 1 \otimes \alpha(x))\\
& = \alpha^2(x) \otimes 1 \otimes 1 + 1 \otimes (\alpha^2(x) \otimes 1 + 1 \otimes \alpha^2(x))\\
& = \alpha^2(x) \otimes 1 \otimes 1 + 1\otimes\alpha^2(x)\otimes 1 + 1 \otimes 1 \otimes \alpha^2(x)\\
& = (\Delta_L \otimes \alpha)(\alpha(x) \otimes 1 + 1 \otimes \alpha(x))\\
& = ((\Delta_L \otimes \alpha) \circ \Delta_L)(x).
\end{split}
\]
This shows that $\Delta_L$ is Hom-coassociative, as desired.
\end{proof}

\section{The Hom-affine plane as an $\bM$-comodule Hom-algebra}
\label{sec:comod}

The main purpose of this section is to generalize to the Hom-algebra setting the following fact about the affine plane: Let $M(2)$ denote the polynomial algebra $\bk[a,b,c,d]$.  It is a bialgebra when equipped with the comultiplication $\Delta \colon M(2) \to M(2)^{\otimes 2}$, defined exactly as in \eqref{Delta}.  Then there is an $M(2)$-comodule algebra structure on the affine plane $A_2 = \bk[x,y]$.  In matrix notation, the $M(2)$-comodule structure map on $A_2$ is given by
\[
\Delta_A \begin{pmatrix}x\\ y\end{pmatrix} = \begin{pmatrix}\Delta_A(x)\\ \Delta_A(y)\end{pmatrix} = \begin{pmatrix}a & b \\ c & d\end{pmatrix} \otimes \begin{pmatrix}x\\ y\end{pmatrix} = \begin{pmatrix}a \otimes x + b \otimes y \\ c \otimes x + d \otimes y\end{pmatrix}.
\]
The reader can consult, e.g., \cite[Theorem III.7.3]{kassel}, for the proof of the above statement.

To extend the above statement to the Hom-algebra setting, we first need to define the concept of \emph{comodule Hom-algebra}.

By a \textbf{multiplicative Hom-coassociative coalgebra}, we mean a triple $(C,\Delta,\alpha)$ in which:
\begin{enumerate}
\item
$(C,\alpha)$ is a Hom-module;
\item
The comultiplication $\Delta \colon C \to C^{\otimes 2}$ is Hom-coassociative, i.e., $(\Delta \otimes \alpha) \circ \Delta = (\alpha \otimes \Delta) \circ \Delta$;
\item
The linear map $\alpha$ is comultiplicative, i.e., $\Delta \circ \alpha = \alpha^{\otimes 2} \circ \Delta$.
\end{enumerate}
Note that our definition of a multiplicative Hom-coassociative coalgebra is slightly different from that of \cite[Definition 1.4]{ms2}.  If $(B,\mu,\Delta,\alpha)$ is a Hom-bialgebra, then clearly $(B,\Delta,\alpha)$ is a multiplicative Hom-coassociative coalgebra.

Let $(C,\Delta_C,\alpha_C)$ be a multiplicative Hom-coassociative coalgebra, and let $(M,\alpha_M)$ be a Hom-module.  A \textbf{$C$-comodule} structure on $M$ consists of a linear map $\Delta_M \colon M \to C \otimes M$ such that $(\Delta_C \otimes \alpha_M) \circ \Delta_M = (\alpha_C \otimes \Delta_M) \circ \Delta_M$.

\begin{definition}
Let $(H,\mu_H,\Delta_H,\alpha_H)$ be a Hom-bialgebra, and let $(A,\mu_A,\alpha_A)$ be a multiplicative Hom-associative algebra.  Then an \textbf{$H$-comodule Hom-algebra} structure on $A$ consists of an $H$-comodule structure $\Delta_A \colon A \to H \otimes A$ on $A$ such that $\Delta_A$ is a morphism of Hom-associative algebras.  We call $\Delta_A$ the \textbf{structure map} of the $H$-comodule Hom-algebra $A$.
\end{definition}

The above definition of an $H$-comodule Hom-algebra coincides with the usual definition of an $H$-comodule algebra when $H$ is a bialgebra and $A$ is an associative algebra (i.e., when $\alpha_H = Id_H$ and $\alpha_A = Id_A$).  The following result shows that comodule algebras deform into comodule Hom-algebras via endomorphisms.  This gives a large class of examples of comodule Hom-algebras.

\begin{theorem}
\label{prop:comod}
Let $(H,\mu_H,\Delta_H)$ be a bialgebra, and let $(A,\mu_A)$ be an $H$-comodule algebra with structure map $\rho \colon A \to H \otimes A$.  Let $\alpha_H \colon H \to H$ be a bialgebra morphism, and let $\alpha_A \colon A \to A$ be an algebra morphism such that \begin{equation}
\label{eq:rhoalpha}
\rho \circ \alpha_A = (\alpha_H \otimes \alpha_A) \circ \rho.
\end{equation}
Define the map $\rho_\alpha = \rho \circ \alpha_A \colon A \to H \otimes A$.  Then:
\begin{enumerate}
\item
$(H,\mu_{\alpha,H} = \alpha_H \circ \mu_H,\Delta_\alpha = \Delta_H \circ \alpha_H, \alpha_H)$ is a Hom-bialgebra, and $(A,\mu_{\alpha,A} = \alpha_A \circ \mu_A,\alpha_A)$ is a multiplicative Hom-associative algebra.
\item
$\rho_\alpha$ is the structure map of an $H$-comodule Hom-algebra structure on $A$, where $H$ and $A$ are given the Hom-bialgebra and multiplicative Hom-associative algebra structures, respectively, of the previous statement.
\end{enumerate}
\end{theorem}

\begin{proof}
The first assertion is from Example ~\ref{ex:bialgdeform}.  To prove the second assertion, we need to show that:
\begin{enumerate}
\item
$\rho_\alpha$ gives the Hom-module $(A,\alpha_A)$ the structure of an $H$-comodule, where $H$ denotes the multiplicative Hom-coassociative coalgebra $(H,\Delta_\alpha = \Delta_H \circ \alpha_H,\alpha_H)$;
\item
$\rho_\alpha$ is a morphism of Hom-associative algebras.
\end{enumerate}
First, $\rho_\alpha$ gives $(A,\alpha_A)$ the structure of an $H$-comodule if and only if
\begin{equation}
\label{eq:rhoalphacomod}
(\alpha_H \otimes \rho_\alpha) \circ \rho_\alpha = (\Delta_\alpha \otimes \alpha_A) \circ \rho_\alpha.
\end{equation}
The following commutative diagram shows that \eqref{eq:rhoalphacomod} is true:
\[
\SelectTips{cm}{10}
\xymatrix{
A \ar[rr]_-{\alpha_A} \ar[dd]_-{\rho_\alpha}
\ar@{<} `u[rrrr]  `_d[rrrr] ^-{\rho_\alpha} [rrrr]
& & A \ar[rr]_-{\rho} \ar[ddll]_-{\rho} \ar[d]^-{\alpha_A} & & H \otimes A \ar[d]_-{\alpha_H \otimes \alpha_A}
\ar@/^2pc/[dd]^-{\alpha_H \otimes \rho_\alpha}\\
 & & A \ar[rr]^-{\rho} \ar[d]^-{\rho} & & H \otimes A \ar[d]_-{Id_H \otimes \rho}\\
H \otimes A \ar[rr]^-{\alpha_H \otimes \alpha_A} \ar@{<} `d[rrrr]  `[rrrr] _-{\Delta_\alpha \otimes \alpha_A} [rrrr]
& & H \otimes A \ar[rr]^-{\Delta_H \otimes Id_A} & & H^{\otimes 2} \otimes A.
}
\]
The lower triangle and the upper-right square are commutative by \eqref{eq:rhoalpha}.  The lower-right square is commutative because $A$ is an $H$-comodule (in the usual sense) by assumption.

That $\rho_\alpha$ is a morphism of Hom-associative algebras follows from the following commutative diagram, where $\nu = (\mu_H \otimes \mu_A) \circ (Id \otimes \tau \otimes Id)$ and $\eta = (\mu_{\alpha,H} \otimes \mu_{\alpha,A}) \circ (Id \otimes \tau \otimes Id)$:
\[
\SelectTips{cm}{10}
\xymatrix{
A \otimes A \ar[rr]_-{\alpha_A^{\otimes 2}} \ar[d]^-{\mu_A} \ar@/_2pc/[dd]_-{\mu_{\alpha,A}}
\ar `u[rrrr] `_d[rrrr] ^-{\rho_\alpha^{\otimes 2}} [rrrr]
& & A \otimes A \ar[rr]_-{\rho^{\otimes 2}} \ar[d]^-{\mu_A} & & H \otimes A \otimes H \otimes A \ar[d]_-{\nu} \ar@/^2pc/[dd]^-{\eta}\\
A \ar[rr]^-{\alpha_A} \ar[d]^-{\alpha_A} & & A \ar[rr]^-{\rho} \ar[d]^-{\alpha_A} & & H \otimes A \ar[d]_-{\alpha_H \otimes \alpha_A}\\
A \ar[rr]^-{\alpha_A} \ar@{<} `d[rrrr]  `[rrrr] _-{\rho_\alpha} [rrrr]
& & A \ar[rr]^-{\rho} & & H \otimes A.
}
\]
The upper-left square is commutative because $\alpha_A$ is a morphism of algebras.  The lower-right square is commutative by \eqref{eq:rhoalpha}.  The upper-right square is commutative because $\rho \colon A \to H \otimes A$ is a morphism of algebras.
\end{proof}

\subsection{The Hom-affine plane $\bA$}

The usual affine plane $A_2$ is the polynomial algebra $\bk[x,y]$, i.e., the free commutative algebra on two generators.  Thus, it makes sense to define the \textbf{Hom-affine plane} as $\bA = F(\bk\langle x,y \rangle)$, where $x$ and $y$ are two independent variables and $F$ is the free multiplicative Hom-associative algebra functor (Corollary ~\ref{cor:F}).

We are now ready to show that $\bA$ is an $\bM$-comodule Hom-algebra.

\begin{theorem}
\label{thm:comodule}
There is an $\bM$-comodule Hom-algebra structure on the Hom-affine plane $\bA$ given by
\begin{equation}
\label{DeltaA}
\Delta_A \begin{pmatrix}x\\ y\end{pmatrix} = \begin{pmatrix}\Delta_A(x)\\ \Delta_A(y)\end{pmatrix} = \begin{pmatrix}a & b \\ c & d\end{pmatrix} \otimes \begin{pmatrix}x\\ y\end{pmatrix} = \begin{pmatrix}a \otimes x + b \otimes y \\ c \otimes x + d \otimes y\end{pmatrix}
\end{equation}
on the generators $x,y \in \bA$.
\end{theorem}

\begin{proof}
By Corollary ~\ref{cor:F}, we have an adjunction
\[
\HA(\bA, \bM\otimes\bA) \cong \Mod(\bk\langle x,y \rangle,\bM\otimes\bA).
\]
Thus, to define a Hom-associative algebra morphism $\Delta_A \colon \bA \to \bM \otimes \bA$, it suffices to define a linear map $\Delta_A \colon \bk\langle x,y \rangle \to \bM\otimes\bA$.  We define such a linear map as in \eqref{DeltaA}.

It remains to show that $\Delta_A$ gives $\bA$ the structure of an $\bM$-comodule.  In other words, we need to check the equality
\[
(\Delta_M \otimes \alpha_A) \circ \Delta_A = (\alpha_M \otimes \Delta_A) \circ \Delta_A.
\]
Here $\Delta_M$ and $\alpha_M$ denote the comultiplication $\Delta$ (defined in \eqref{Delta}) and $\alpha$ in $\bM$, respectively.  Since $\Delta_A$ and $\Delta_M$ are morphisms of Hom-associative algebras and since $\alpha_A$ and $\alpha_M$ are both multiplicative, it suffices to check this equality on the generators $x$ and $y$ in $\bA$.    We will use the identification
\begin{equation}
\label{MMA}
\bM \otimes \bM \otimes \bA  \cong F(\bk\langle a',a'',b',b'',c',c'',d',d'',x,y\rangle)
\end{equation}
of multiplicative Hom-associative algebras provided by Lemma ~\ref{lem1:F}.  Then, on the one hand, we have
\[
\begin{split}
((\Delta_M \otimes \alpha_A) \circ \Delta_A)\begin{pmatrix}x\\ y\end{pmatrix}
& = \left[\begin{pmatrix}a' & b'\\ c' & d'\end{pmatrix}
\begin{pmatrix}a'' & b''\\ c'' & d''\end{pmatrix}\right]\begin{pmatrix}\alpha(x) \\ \alpha(y)\end{pmatrix}\\
& = \begin{pmatrix}(a'a'')\alpha(x) + (b'c'')\alpha(x) + (a'b'')\alpha(y) + (b'd'')\alpha(y)\\ (c'a'')\alpha(x) + (d'c'')\alpha(x) + (c'b'')\alpha(y) + (d'd'')\alpha(y)\end{pmatrix}.
\end{split}
\]
On the other hand, we have
\[
\begin{split}
((\alpha_M \otimes \Delta_A) \circ \Delta_A)\begin{pmatrix}x\\ y\end{pmatrix}
& = \begin{pmatrix}\alpha(a') & \alpha(b')\\ \alpha(c') & \alpha(d')\end{pmatrix}\left[\begin{pmatrix}a'' & b''\\ c'' & d''\end{pmatrix}\begin{pmatrix}x\\ y\end{pmatrix}\right]\\
& = \begin{pmatrix}\alpha(a')(a''x) + \alpha(b')(c''x) + \alpha(a')(b''y) + \alpha(b')(d''y)\\ \alpha(c')(a''x) + \alpha(d')(c''x) + \alpha(c')(b''y) + \alpha(d')(d''y)\end{pmatrix}.
\end{split}
\]
Using the Hom-associativity in the Hom-associative algebra in \eqref{MMA}, we conclude that $(\Delta_M \otimes \alpha_A) \circ \Delta_A$ and $(\alpha_M \otimes \Delta_A) \circ \Delta_A$ are equal.  This shows that $\Delta_A$ as in \eqref{DeltaA} gives the Hom-affine plane $\bA$ the structure of an $\bM$-comodule Hom-algebra.
\end{proof}


\end{document}